\documentstyle[12pt]{article}

\newcommand{\Conj}{\mathop{\bigwedge\kern -6pt\bigwedge}}

\newcommand{\forces}{\mid\vdash}

\newtheorem{theorem}{Theorem}
\newtheorem{lemma}[theorem]{Lemma}

\newcommand{\res}{\upharpoonright}

\catcode`\@=11

\font\teneuf=eufm10  scaled 1200
\font\seveneuf=eufm7 scaled 1200
\font\fiveeuf=eufm5  scaled 1200
\font\tenmsx=msxm10  scaled 1200
\font\sevenmsx=msxm7 scaled 1200
\font\fivemsx=msxm5  scaled 1200
\font\tenmsy=msym10  scaled 1200
\font\sevenmsy=msym7 scaled 1200
\font\fivemsy=msym5  scaled 1200

\newfam\euffam
\newfam\msxfam
\newfam\msyfam

\textfont\euffam=\teneuf \scriptfont\euffam=\seveneuf
 \scriptscriptfont\euffam=\fiveeuf
\textfont\msxfam=\tenmsx  \scriptfont\msxfam=\sevenmsx
  \scriptscriptfont\msxfam=\fivemsx
\textfont\msyfam=\tenmsy  \scriptfont\msyfam=\sevenmsy
  \scriptscriptfont\msyfam=\fivemsy

\def\frak{\ifmmode\let\next\frak@\else
 \def\next{\errmessage{Use \string\frak\space only in math mode}}\fi\next}
\def\goth{\ifmmode\let\next\frak@\else
 \def\next{\errmessage{Use \string\goth\space only in math mode}}\fi\next}
\def\frak@#1{{\frak@@{#1}}}
\def\frak@@#1{\fam\euffam#1}

\def\Bbb{\ifmmode\let\next\Bbb@\else
 \def\next{\errmessage{Use \string\Bbb\space only in math mode}}\fi\next}
\def\Bbb@#1{{\Bbb@@{#1}}}
\def\Bbb@@#1{\fam\msyfam#1}

\def\hexnumber@#1{\ifcase#1 0\or1\or2\or3\or4\or5\or6\or7\or8\or9\or
        A\or B\or C\or D\or E\or F\fi }

\edef\msx@{\hexnumber@\msxfam}
\edef\msy@{\hexnumber@\msyfam}

\mathchardef\boxdot="2\msx@00
\mathchardef\boxplus="2\msx@01
\mathchardef\boxtimes="2\msx@02
\mathchardef\square="0\msx@03
\mathchardef\blacksquare="0\msx@04
\mathchardef\centerdot="2\msx@05
\mathchardef\lozenge="0\msx@06
\mathchardef\blacklozenge="0\msx@07
\mathchardef\circlearrowright="3\msx@08
\mathchardef\circlearrowleft="3\msx@09
\mathchardef\rightleftharpoons="3\msx@0A
\mathchardef\leftrightharpoons="3\msx@0B
\mathchardef\boxminus="2\msx@0C
\mathchardef\Vdash="3\msx@0D
\mathchardef\Vvdash="3\msx@0E
\mathchardef\vDash="3\msx@0F
\mathchardef\twoheadrightarrow="3\msx@10
\mathchardef\twoheadleftarrow="3\msx@11
\mathchardef\leftleftarrows="3\msx@12
\mathchardef\rightrightarrows="3\msx@13
\mathchardef\upuparrows="3\msx@14
\mathchardef\downdownarrows="3\msx@15
\mathchardef\upharpoonright="3\msx@16
\mathchardef\downharpoonright="3\msx@17
\mathchardef\upharpoonleft="3\msx@18
\mathchardef\downharpoonleft="3\msx@19
\mathchardef\rightarrowtail="3\msx@1A
\mathchardef\leftarrowtail="3\msx@1B
\mathchardef\leftrightarrows="3\msx@1C
\mathchardef\rightleftarrows="3\msx@1D
\mathchardef\Lsh="3\msx@1E
\mathchardef\Rsh="3\msx@1F
\mathchardef\rightsquigarrow="3\msx@20
\mathchardef\leftrightsquigarrow="3\msx@21
\mathchardef\looparrowleft="3\msx@22
\mathchardef\looparrowright="3\msx@23
\mathchardef\circeq="3\msx@24
\mathchardef\succsim="3\msx@25
\mathchardef\gtrsim="3\msx@26
\mathchardef\gtrapprox="3\msx@27
\mathchardef\multimap="3\msx@28
\mathchardef\therefore="3\msx@29
\mathchardef\because="3\msx@2A
\mathchardef\doteqdot="3\msx@2B
\mathchardef\triangleq="3\msx@2C
\mathchardef\precsim="3\msx@2D
\mathchardef\lesssim="3\msx@2E
\mathchardef\lessapprox="3\msx@2F
\mathchardef\eqslantless="3\msx@30
\mathchardef\eqslantgtr="3\msx@31
\mathchardef\curlyeqprec="3\msx@32
\mathchardef\curlyeqsucc="3\msx@33
\mathchardef\preccurlyeq="3\msx@34
\mathchardef\leqq="3\msx@35
\mathchardef\leqslant="3\msx@36
\mathchardef\lessgtr="3\msx@37
\mathchardef\backprime="0\msx@38
\mathchardef\risingdotseq="3\msx@3A
\mathchardef\fallingdotseq="3\msx@3B
\mathchardef\succcurlyeq="3\msx@3C
\mathchardef\geqq="3\msx@3D
\mathchardef\geqslant="3\msx@3E
\mathchardef\gtrless="3\msx@3F
\mathchardef\sqsubset="3\msx@40
\mathchardef\sqsupset="3\msx@41
\mathchardef\vartriangleright="3\msx@42
\mathchardef\vartriangleleft="3\msx@43
\mathchardef\trianglerighteq="3\msx@44
\mathchardef\trianglelefteq="3\msx@45
\mathchardef\bigstar="0\msx@46
\mathchardef\between="3\msx@47
\mathchardef\blacktriangledown="0\msx@48
\mathchardef\blacktriangleright="3\msx@49
\mathchardef\blacktriangleleft="3\msx@4A
\mathchardef\vartriangle="0\msx@4D
\mathchardef\blacktriangle="0\msx@4E
\mathchardef\triangledown="0\msx@4F
\mathchardef\eqcirc="3\msx@50
\mathchardef\lesseqgtr="3\msx@51
\mathchardef\gtreqless="3\msx@52
\mathchardef\lesseqqgtr="3\msx@53
\mathchardef\gtreqqless="3\msx@54
\mathchardef\Rrightarrow="3\msx@56
\mathchardef\Lleftarrow="3\msx@57
\mathchardef\veebar="2\msx@59
\mathchardef\barwedge="2\msx@5A
\mathchardef\doublebarwedge="2\msx@5B
\mathchardef\angle="0\msx@5C
\mathchardef\measuredangle="0\msx@5D
\mathchardef\sphericalangle="0\msx@5E
\mathchardef\varpropto="3\msx@5F
\mathchardef\smallsmile="3\msx@60
\mathchardef\smallfrown="3\msx@61
\mathchardef\Subset="3\msx@62
\mathchardef\Supset="3\msx@63
\mathchardef\Cup="2\msx@64
\mathchardef\Cap="2\msx@65
\mathchardef\curlywedge="2\msx@66
\mathchardef\curlyvee="2\msx@67
\mathchardef\leftthreetimes="2\msx@68
\mathchardef\rightthreetimes="2\msx@69
\mathchardef\subseteqq="3\msx@6A
\mathchardef\supseteqq="3\msx@6B
\mathchardef\bumpeq="3\msx@6C
\mathchardef\Bumpeq="3\msx@6D
\mathchardef\lll="3\msx@6E
\mathchardef\ggg="3\msx@6F
\mathchardef\circledS="0\msx@73
\mathchardef\pitchfork="3\msx@74
\mathchardef\dotplus="2\msx@75
\mathchardef\backsim="3\msx@76
\mathchardef\backsimeq="3\msx@77
\mathchardef\complement="0\msx@7B
\mathchardef\intercal="2\msx@7C
\mathchardef\circledcirc="2\msx@7D
\mathchardef\circledast="2\msx@7E
\mathchardef\circleddash="2\msx@7F
\def\ulcorner{\delimiter"4\msx@70\msx@70 }
\def\urcorner{\delimiter"5\msx@71\msx@71 }
\def\llcorner{\delimiter"4\msx@78\msx@78 }
\def\lrcorner{\delimiter"5\msx@79\msx@79 }
\def\yen{\mathhexbox\msx@55 }
\def\checkmark{\mathhexbox\msx@58 }
\def\circledR{\mathhexbox\msx@72 }
\def\maltese{\mathhexbox\msx@7A }
\mathchardef\lvertneqq="3\msy@00
\mathchardef\gvertneqq="3\msy@01
\mathchardef\nleq="3\msy@02
\mathchardef\ngeq="3\msy@03
\mathchardef\nless="3\msy@04
\mathchardef\ngtr="3\msy@05
\mathchardef\nprec="3\msy@06
\mathchardef\nsucc="3\msy@07
\mathchardef\lneqq="3\msy@08
\mathchardef\gneqq="3\msy@09
\mathchardef\nleqslant="3\msy@0A
\mathchardef\ngeqslant="3\msy@0B
\mathchardef\lneq="3\msy@0C
\mathchardef\gneq="3\msy@0D
\mathchardef\npreceq="3\msy@0E
\mathchardef\nsucceq="3\msy@0F
\mathchardef\precnsim="3\msy@10
\mathchardef\succnsim="3\msy@11
\mathchardef\lnsim="3\msy@12
\mathchardef\gnsim="3\msy@13
\mathchardef\nleqq="3\msy@14
\mathchardef\ngeqq="3\msy@15
\mathchardef\precneqq="3\msy@16
\mathchardef\succneqq="3\msy@17
\mathchardef\precnapprox="3\msy@18
\mathchardef\succnapprox="3\msy@19
\mathchardef\lnapprox="3\msy@1A
\mathchardef\gnapprox="3\msy@1B
\mathchardef\nsim="3\msy@1C

\mathchardef\ncong="3\msy@1D
\mathchardef\varsubsetneq="3\msy@20
\mathchardef\varsupsetneq="3\msy@21
\mathchardef\nsubseteqq="3\msy@22
\mathchardef\nsupseteqq="3\msy@23
\mathchardef\subsetneqq="3\msy@24
\mathchardef\supsetneqq="3\msy@25
\mathchardef\varsubsetneqq="3\msy@26
\mathchardef\varsupsetneqq="3\msy@27
\mathchardef\subsetneq="3\msy@28
\mathchardef\supsetneq="3\msy@29
\mathchardef\nsubseteq="3\msy@2A
\mathchardef\nsupseteq="3\msy@2B
\mathchardef\nparallel="3\msy@2C
\mathchardef\nmid="3\msy@2D
\mathchardef\nshortmid="3\msy@2E
\mathchardef\nshortparallel="3\msy@2F
\mathchardef\nvdash="3\msy@30
\mathchardef\nVdash="3\msy@31
\mathchardef\nvDash="3\msy@32
\mathchardef\nVDash="3\msy@33
\mathchardef\ntrianglerighteq="3\msy@34
\mathchardef\ntrianglelefteq="3\msy@35
\mathchardef\ntriangleleft="3\msy@36
\mathchardef\ntriangleright="3\msy@37
\mathchardef\nleftarrow="3\msy@38
\mathchardef\nrightarrow="3\msy@39
\mathchardef\nLeftarrow="3\msy@3A
\mathchardef\nRightarrow="3\msy@3B
\mathchardef\nLeftrightarrow="3\msy@3C
\mathchardef\nleftrightarrow="3\msy@3D
\mathchardef\divideontimes="2\msy@3E
\mathchardef\varnothing="0\msy@3F
\mathchardef\nexists="0\msy@40
\mathchardef\mho="0\msy@66
\mathchardef\eth="0\msy@67
\mathchardef\eqsim="3\msy@68
\mathchardef\beth="0\msy@69
\mathchardef\gimel="0\msy@6A
\mathchardef\daleth="0\msy@6B
\mathchardef\lessdot="3\msy@6C
\mathchardef\gtrdot="3\msy@6D
\mathchardef\ltimes="2\msy@6E
\mathchardef\rtimes="2\msy@6F
\mathchardef\shortmid="3\msy@70
\mathchardef\shortparallel="3\msy@71
\mathchardef\smallsetminus="2\msy@72
\mathchardef\thicksim="3\msy@73
\mathchardef\thickapprox="3\msy@74
\mathchardef\approxeq="3\msy@75
\mathchardef\succapprox="3\msy@76
\mathchardef\precapprox="3\msy@77
\mathchardef\curvearrowleft="3\msy@78
\mathchardef\curvearrowright="3\msy@79
\mathchardef\digamma="0\msy@7A
\mathchardef\varkappa="0\msy@7B
\mathchardef\hslash="0\msy@7D
\mathchardef\hbar="0\msy@7E
\mathchardef\backepsilon="3\msy@7F

\catcode`\@=12

\begin{document}

  \begin{center}
     {\large A MODEL IN WHICH THERE ARE JECH--KUNEN TREES\\ 
             BUT THERE ARE NO KUREPA TREES}
  \footnote{1980 Mathematics Subject Classification (1985 Revision).
	    Primary 03E35.}
  \end{center}

  \begin{center}
     Saharon Shelah \footnote{The research of the first author
was partially supported by the
Basic Research Fund, Israeli Acad. of Science Publ. nu. 466.} and Renling Jin
  \end{center}

  \bigskip

  \begin{quote}

    \centerline{Abstract}

    \small

By an $\omega_{1}$--tree we mean a tree of power $\omega_{1}$ and height
$\omega_{1}$. We call an $\omega_{1}$--tree a Jech--Kunen tree if it has
$\kappa$--many branches for some $\kappa$ strictly between $\omega_{1}$
and $2^{\omega_{1}}$. In this paper we construct the models of $C\!H$
plus $2^{\omega_{1}}>\omega_{2}$, in which there are Jech--Kunen trees
and there are no Kurepa trees.

  \end{quote}

An partially ordered set, or poset for short, $\langle T,<_{T}\rangle$
is called a tree if for every $t\in T$ the set $\{s\in T:s<_{T}t\}$ is
well--ordered under $<_{T}$. The order type of that set is called the
height of $t$ in $T$, denoted by $ht_{T}(t)$. We will not distinguish
a tree from its base set. For every ordinal $\alpha$, let $T_{\alpha}$,
the $\alpha$--th level of $T$, $=\{t\in T:ht_{T}(t)=\alpha\}$ and
$T\!\res\!\alpha =\bigcup_{\beta<\alpha}T_{\beta}$. Let $ht(T)$, the height
of $T$, is the smallest ordinal $\alpha$ such that $T_{\alpha}=\emptyset$.
By a branch of $T$ we mean a linearly ordered subset of $T$ which 
intersects every nonempty level of $T$. Let ${\cal B}(T)$ be the set of
all branches of $T$. $T'$ is called a subtree of $T$ if $T'\subseteq T$
and $<_{T'}=<_{T}\bigcap T'\times T'$ ($T'$ inherits the order of $T$).

$T$ is called an $\omega_{1}$--tree if $|T|=\omega_{1}$ and $ht(T)=\omega_{1}$.
An $\omega_{1}$--tree $T$ is called a Kurepa tree if 
$|{\cal B}(T)|>\omega_{1}$ and for every $\alpha\in\omega_{1}$,
$|T_{\alpha}|<\omega_{1}$.
An $\omega_{1}$--tree is called a Jech--Kunen tree if $\omega_{1}<
|{\cal B}(T)|<2^{\omega_{1}}$.

\medskip

The independence of the existence of Kurepa trees was proved by J. H. Silver
(see [K2, \S 3 of Chapter VIII]). T. Jech in [Je1] constructed by forcing a
model of $C\!H$ plus $2^{\omega_{1}}>\omega_{2}$, in which there is a
Jech--Kunen tree. In fact, it is a Kurepa tree with fewer than 
$2^{\omega_{1}}$--many branches. The independence of the existence of
Jech--Kunen trees under $C\!H$ plus $2^{\omega_{1}}>\omega_{2}$ was given
by K. Kunen [K1]. In his paper he gave an equivalent form of Jech--Kunen trees
in terms of compact Hausdorff spaces. The detailed proof can be found in
[Ju, Theorem 4.8].

In both Silver and Kunen's proofs, the existence of a strongly inaccessible
cardinal was assumed (the assumption is also necessary). The technique
they used to kill all Kurepa trees or Jech--Kunen trees is to show that
if an $\omega_{1}$--tree $T$ has a new branch in an $\omega_{1}$--closed
forcing extension, then $T$ must have a subtree which is isomorphic
to $\langle 2^{<\omega_{1}},\subseteq\rangle$, a complete binary tree of
height $\omega_{1}$. So in Kunen's model not only all Jech--Kunen trees
are killed, but also all Kurepa trees are killed.

R. Jin in [Ji1] started discussing the differences between Kurepa trees
and Jech--Kunen trees. He showed that it is independent of $C\!H$ plus
$2^{\omega_{1}}>\omega_{2}$ that there exists a Kurepa tree which has
no Jech--Kunen subtrees. He also showed that it is independent of
$C\!H$ plus $2^{\omega_{1}}>\omega_{2}$ that there exists a Jech--Kunen
tree which has no Kurepa subtrees. In his proofs some strongly inaccessible
cardinals were assumed and later, Kunen eliminated the large cardinal
assumption for one of the proofs.

In [Ji2] Jin proved that assuming the existence of two inaccessible
cardinals, it is consistent with $C\!H$ plus $2^{\omega_{1}}>\omega_{2}$
that there exist Kurepa trees and there are no Jech--Kunen trees.

The problem whether $C\!H$ plus $2^{\omega_{1}}>\omega_{2}$ is consistent
with that there exist Jech--Kunen trees and there are no Kurepa trees,
was posed in [Ji2]. We will answer the question in this paper by assuming
naturally the existence of a strongly inaccessible cardinal.

Before proving our results we need more notations and definitions.

\medskip

A tree $T$ is called normal if,

(1) every $t\in T$ has at least two
immediate successors, 

(2) for every $t\in T$ and an ordinal $\alpha$ such that
$ht_{T}(t)<\alpha<ht(T)$, there exists $t'\in T_{\alpha}$ such that
$t<_{T}t'$. 

A tree $C=\{c_{s}:s\in 2^{<\omega}\}$ is called a Cantor
tree if the map $s\mapsto c_{s}$ is an isomorphism from
$\langle 2^{<\omega},\subseteq\rangle$ to $C$. For convenience we assume,
from now on, that every tree considered in this paper is a subtree of
$\langle 2^{<\omega_{1}},\subseteq\rangle$ with the unique root $\emptyset$.
By that way we can define the least upper bound of an increasing sequence
in a tree by taking its union. Let $\lim (\omega_{1})$ be the set of all
limit ordinals in $\omega_{1}$. Let $T$ be a tree
and $\delta\in\lim (\omega_{1})$. A subtree $C$ of $T$ is called cofinal
in $T\!\res\!\delta$ if for every $B\in {\cal B}(C)$, the set $\{ht_{T}(t):
t\in B\}$ is cofinal in $\delta$. $T$ is called complete at level $\delta$
if for every $B\in {\cal B}(T\!\res\!\delta)$, $\bigcup B\in T_{\delta}$.
$T$ is called properly pruned at level $\delta$ if for every Cantor
subtree $C=\{c_{s}:s\in 2^{<\omega}\}$ of $T$
which is cofinal in $T\!\res\!\delta$, there exist $f, g\in 
2^{\omega}$ such that $\bigcup_{n\in\omega}c_{f\res n}\in T_{\delta}$
and $\bigcup_{n\in\omega}c_{g\res n}\not\in T_{\delta}$

Let $S\subseteq\lim(\omega_{1})$. A tree is called $S$--properly pruned
if for every $\alpha\in\lim(\omega_{1})$, 
$\alpha\not\in S$ implies that $T$ is complete at level
$\alpha$, and $\alpha\in S$ implies that $T$ is properly pruned at level 
$\alpha$.

Let $I$ be an index set and $T$ be a tree. 
For every $F\in T^{I}$, let $supt(F)$, the support
of $F$, be the set $\{i\in I:F(i)\neq\emptyset\}$. Let $F,G\in T^{I}$.
Define $F\preccurlyeq G$ iff for every $i\in I$, $F(i)\leq G(i)$. We call $F\in
T^{I}$ uniform at $\delta$ for some $\delta\in\omega_{1}$ if
for every $i\in supt(F)$, $ht_{T}(F(i))=\delta$. Let $C=\{F_{s}:
s\in 2^{<\omega}\}\subseteq T^{I}$ be a Cantor tree (under $\preccurlyeq$).
$C$ is called uniformly cofinal in $(T\!\res\!\delta)^{I}$ for some $\delta\in
\omega_{1}$ if for every $s\in 2^{<\omega}$, 
there is a $\delta_{s}<\delta$ such that $F_{s}$ is uniform at $\delta_{s}$
and for every $i\in \bigcup_{s\in 2^{<\omega}}supt(F_{s})$, the subtree
$\{F_{s}(i):s\in 2^{<\omega}\}$ of $T$ is cofinal in $T\!\res\!\delta$.
We use $\perp$ for the word ``incompatible''. 
For example, for any $s,t\in 2^{<\omega}$, $s\perp t$ means $s\bigcup t$
is not a function. For any $F, G\in T^{I}$, we call that $F$ and $G$
are completely incompatible if
for any $i\in supt(F)$ and any $j\in supt(G)$, $F(i)\perp G(j)$ 
($F(i)$ and $G(j)$ have no common upper bound in $T$).
Now $C$ is called separated if for any $s,s'\in 2^{<\omega}$, $s\perp s'$
implies that $F_{s}$ and $F_{s'}$ are completely incompatible.

Let $T$ be a tree and $\delta\in\lim (\omega_{1})$. We call that $T$ is
properly pruned in countable products at level $\delta$ if for every
Cantor tree $C=\{F_{s}:s\in 2^{<\omega}\}\subseteq T^{I}$, which is
separated and uniformly cofinal in $(T\!\res\!\delta)^{\omega}$, there exist
$f,g\in 2^{\omega}$ such that for every $i\in\bigcup_{n\in\omega}
supt(F_{f\res n})$, $\bigcup_{n\in\omega}F_{f\res n}(i)\in T_{\delta}$ and
for every $i\in\bigcup_{n\in\omega}supt(F_{g\res n})$, $\bigcup_{n\in\omega}
F_{g\res n}(i)\not\in T_{\delta}$.

Let $S\subseteq\lim(\omega_{1})$. A tree is called $S$--properly pruned
in countable products if for every $\alpha\in\lim (\omega_{1})$,
$\alpha\not\in S$ implies that $T$ is complete at level $\alpha$, 
and $\alpha\in S$ implies that $T$ is properly
pruned in countable products at level $\alpha$.

\begin{lemma}

Let $T$ be a tree and $I$ be an index set. For any Cantor tree $C=
\{F_{s}:s\in 2^{<\omega}\}\subseteq T^{I}$, if $C$ is separated, then
for any $f,g\in 2^{\omega}$, $f\neq g$ implies that
$\langle\bigcup_{n\in\omega}F_{f\res n}(i)\rangle_{i\in I}$ and
$\langle\bigcup_{n\in\omega}F_{g\res n}(i)\rangle_{i\in I}$ are 
completely incompatible.

\end{lemma}

\noindent {\bf Proof:}\quad
Let $i\in\bigcup_{n\in\omega}supt(F_{f\res n})$ and $j\in\bigcup_{n\in\omega}
supt(F_{g\res n})$. Let $m\in\omega$ such that $i\in supt(F_{f\res m})$, 
$j\in supt(F_{g\res m})$ and $f\!\res\! m\neq g\!\res\! m$. 
Then $\bigcup_{n\in\omega}F_{f\res n}(i)$ and $\bigcup_{n\in\omega}
F_{g\res n}(j)$ are compatible implies that $F_{f\res m}(i)$ and 
$F_{g\res m}(j)$ are compatible, a contradiction. \quad $\Box$

\begin{lemma}

($C\!H$). For any $S\subseteq \lim (\omega_{1})$, there exists a normal
$\omega_{1}$--tree which is $S$--properly pruned in countable products.

\end{lemma}

\noindent {\bf Proof:}\quad
We construct $T_{\delta}\subseteq 2^{\delta}$ recursively on
$\delta<\omega_{1}$ and $T=\bigcup_{\delta<\omega_{1}}T_{\delta}$
will be the tree we want.  
 
\medskip

Case 1. \quad $\delta=\beta+1$ for some $\beta\in\omega_{1}$.

Let $T_{\delta}=\{t\hat{\;}\langle l\rangle:t\in T_{\beta},l=0,1\}$.

\medskip

Case 2. \quad $\delta\in\lim (\omega_{1})\smallsetminus S$.

Let $T_{\delta}=\{\bigcup B:B\in {\cal B}(T\!\res\!\delta)\}$.

\medskip

Case 3. \quad $\delta\in S$.

Let $\cal C$ be the set of all Cantor trees which are separated and
uniformly cofinal in $(T\!\res\!\delta)^{\omega}$. By $C\!H$ we have that
$|{\cal C}|\leq (\omega_{1}^{\omega})^{\omega}=\omega_{1}$. Let
${\cal C}=\{C^{\alpha}:\alpha\in\omega_{1}\}$ be an enumeration, where
$C^{\alpha}=\{F^{\alpha}_{s}:s\in 2^{<\omega}\}$. We now want to find
a set $X\subseteq\{\bigcup B:B\in {\cal B}(T\!\res\!\delta)\}$ such that
for every $\alpha\in\omega_{1}$, there are $f,g\in 2^{\omega}$ such that
\[\{\bigcup_{n\in\omega}F^{\alpha}_{f\res n}(i):i\in\omega\}\subseteq
X\bigcup\{\emptyset\}\] and 
\[\{\bigcup_{n\in\omega}F^{\alpha}_{g\res n}(i):
i\in\omega\}\bigcap X=\emptyset.\]
If $X$ is found, we let $T_{\delta}=X$.

We now build $X_{\gamma}$ and $Y_{\gamma}$ recursively such that,

(1) $X_{\gamma}$ and $Y_{\gamma}$ are countable, 

(2) $\gamma<\gamma'<\omega_{1}$ implies that $X_{\gamma}\subseteq
X_{\gamma'}$ and $Y_{\gamma}\subseteq Y_{\gamma'}$,

(3) $X_{\gamma}\bigcap
Y_{\gamma}=\emptyset$ for every $\gamma\in\omega_{1}$, 

(4) for every $\gamma\in\omega_{1}$, there exist $f,g\in 2^{\omega}$ such
that $\{\bigcup_{n\in\omega}F^{\gamma}_{f\res n}(i):i\in\omega\}\subseteq
X_{\gamma+1}$ and $\{\bigcup_{n\in\omega}F^{\gamma}_{g\res n}(i):i\in
\omega\}\subseteq Y_{\gamma+1}$.

Let $X_{0}=Y_{0}=
\emptyset$. Let $X_{\gamma}=\bigcup_{\beta<\gamma}X_{\beta}$ and
$Y_{\gamma}=\bigcup_{\beta<\gamma}Y_{\beta}$ if $\gamma\in\lim (\omega_{1})$.
For $\gamma+1$, since $X_{\gamma}$ and $Y_{\gamma}$ are countable and
$C^{\gamma}$ is separated, 
by Lemma 1, there exist $f,g\in 2^{\omega},f\neq g$ such that
\[(X_{\gamma}\bigcup Y_{\gamma})\bigcap (\{\bigcup_{n\in\omega}F^{\gamma}
_{f\res n}(i):i\in\omega\}\bigcup\{\bigcup_{n\in\omega}F^{\gamma}_{g\res n}(i)
:i\in\omega\})=\emptyset.\] Hence let 
\[X_{\gamma+1}=X_{\gamma}\bigcup
\{\bigcup_{n\in\omega}F^{\gamma}_{f\res n}(i):i\in\bigcup_{n\in\omega}
supt(F^{\gamma}_{f\res n})\}\] and 
\[Y_{\gamma+1}=Y_{\gamma}\bigcup
\{\bigcup_{n\in\omega}F^{\gamma}_{g\res n}(i):i\in\bigcup_{n\in\omega}
supt(F^{\gamma}_{g\res n})\}.\] Then $X=\bigcup_{\gamma\in\omega_{1}}
X_{\gamma}$ is the set we want. \quad $\Box$

\begin{lemma}

Let $S\subseteq\lim (\omega_{1})$. $T$ is $S$--properly pruned in
countable products implies that $T$ is $S$--properly pruned.

\end{lemma}

\noindent {\bf Proof:}\quad
If $C=\{c_{s}:s\in 2^{<\omega}\}\subseteq T$ is a Cantor tree which is cofinal
in $T\!\res\!\delta$ for some $\delta\in S$, then the Cantor tree $D=\{F_{s}:
s\in 2^{<\omega}\}\subseteq T^{\omega}$, where $F_{s}(0)=c_{s}$ 
and $F_{s}(i)=\emptyset$ for every
$i\neq 0$, is separated and uniformly cofinal in $(T\!\res\!\delta)^{\omega}$.
\quad $\Box$

\begin{lemma}

Let $S\subseteq\lim (\omega_{1})$ and $T$ be $S$--properly pruned in
countable products. Let $C=\{F_{s}:s\in 2^{<\omega}\}$ be a separated and
uniformly cofinal Cantor subtree in $(T\!\res\!\delta)^{\omega}$ for some
$\delta\in S$. Then there are uncountably many $f\in 2^{\omega}$ such
that for every $i\in\bigcup_{n\in\omega}supt(F_{f\res n})$,
$\bigcup_{n\in\omega}F_{f\res n}(i)\in T_{\delta}$.

\end{lemma}

\noindent {\bf Proof:}\quad
Suppose that the lemma is not true. Then we can find a Cantor subtree
$C'=\{F'_{s}:s\in 2^{<\omega}\}\subseteq C$ such that for every $f\in
2^{\omega}$, there exists $i\in\omega$, $\bigcup_{n\in\omega}F'_{f\res n}(i)
\not\in T_{\delta}$. Since $C'$ is a subtree of $C$, $C'$ itself
is also separated and uniformly cofinal
in $(T\!\res\!\delta)^{\omega}$.
That contradicts the definition of the $S$--properly
prunedness in countable products. \quad $\Box$

\bigskip

Next we shall use the forcing method to construct desired models.
For the terminology and basic facts of forcing, see [K2] and [Je2].
We always assume the consistency of $Z\!F\!C$ and let $M$ be always a
countable transitive model of $Z\!F\!C$. In the forcing arguments,
we always let $\dot{a}$ be a name of $a$ if $a$ is not in the ground
model. For every element $a$ in the ground model, we will not distinguish 
$a$ from its canonical name.

Let $I,J$ be two sets. Let 
\[Fn(I,J,\omega_{1})=\{p:p\subseteq I\times J
\mbox{ is a function and }|p|<\omega_{1}\}\] 
be a poset ordered by
reverse inclusion. Let $I$ be a subset of a cardinal $\kappa$. Let
\[Lv(I,\omega_{1})=\]
\[\{p:p\subseteq (I\times\omega_{1})\times\kappa
\mbox{ is a function, }|p|<\omega_{1}\mbox{ and }
\forall\langle\alpha,\beta\rangle\in\mbox{dom}(p)
(p(\alpha,\beta)\in\alpha)\}\]
be a poset ordered by reverse inclusion.
Let $T$ be a tree and $I$ be an index set. Let 
\[{\Bbb P}(T,I,\omega_{1})
=\{F:F\in T^{I}\mbox{ and }|supt(F)|<\omega_{1}\}.\]
The order of ${\Bbb P}
(T,I,\omega_{1})$ is defined as the reverse order of $T^{I}$, or
$F\leq_{{\Bbb P}(T,I,\omega_{1})}G$ iff $G\preccurlyeq F$.

\begin{lemma}

Let $T$ be a normal $\omega_{1}$--tree and $I$ be an index set.
For any $p,q\in {\Bbb P}(T,I,\omega_{1})$, there exist $p',q'\in
{\Bbb P}(T,I,\omega_{1})$ such that $p'<p$ and $q'<q$, $p',q'$
are uniform at $\delta$ for some $\delta\in\omega_{1}$, and
$p'$ is completely incompatible with $q'$.

\end{lemma}

\noindent {\bf Proof:}\quad
Let $\alpha\in\omega_{1}$ be large enough so that $p,q\in 
(T\!\res\!\alpha)^{I}$ ($\alpha$ exists 
because $p,q$ both have countable supports).
Let $\delta>\alpha$ be countable 
such that for every $i\in supt(p)$ 
\[|\{t\in T_{\delta}:
p(i)<_{T}t\}|\geq\omega,\]
and for every $j\in supt(q)$ 
\[|\{t\in T_{\delta}:
q(j)<_{T}t\}|\geq\omega.\]
$\delta$ exists because $T$ is normal.
Let \[supt(p)=\{i_{n}:n\in\omega\}\]
and \[supt(q)=\{j_{n}:n\in\omega\}.\]
We now define $p'(i_{n})$ and $q'(j_{n})$ such that
\[p'(i_{n}),q'(j_{n})\in T_{\delta},\] \[p'(i_{n})>p(i_{n}),\; 
q'(j_{n})>q(j_{n}),\;p'(i_{n})\neq q'(j_{n})\] and
\[p'(i_{n}),q'(j_{n})\not\in\{p'(i_{m}),q'(j_{m}):m<n\}.\]
Let $p'(i)=\emptyset$ if $i\not\in supt(p)$ and let $q'(j)=\emptyset$
if $j\not\in supt(q)$. Then $p'$ and $q'$ are the desired elements.
\quad $\Box$

\bigskip

Let ${\Bbb P}$ be a poset and $D\subseteq {\Bbb P}$. $D$ is called dense
in $\Bbb P$ if for every $p\in {\Bbb P}$ there is $d\in D$ such that
$d\leq p$. $D$ is called open in $\Bbb P$ if for every $p\in {\Bbb P}$
and $d\in D$, $p\leq d$ implies that $p\in D$. $\Bbb P$ is called
$\omega_{1}$--Baire if for any countable sequence 
$\langle D_{n}:n\in\omega\rangle$ of dense open subsets of $\Bbb P$,
$\bigcap_{n\in\omega}D_{n}$ is dense in $\Bbb P$.

\begin{lemma}

In $M$ let $\Bbb P$ be a poset which is $\omega_{1}$--Baire. 
Let $G$ be a $\Bbb P$--generic filter over $M$. Then 
$M^{\omega}\bigcap M[G]\subseteq M$.

\end{lemma}

\noindent {\bf Proof:}\quad
Let $h\in M[G]$ be a function from $\omega$ to $A$, where $A\in M$.
We work in $M$ and let
$p\in G$ such that 
\[p\forces(\dot{h}\mbox{ is a function from }
\omega\mbox{ to }A).\] 
For every $n\in\omega$, let \[D_{n}=\{q\in {\Bbb P}:
q\perp p\mbox{ or }\exists a\in A(q\forces\dot{h}(n)=a)\}.\] Then
$D_{n}$ is dense open in $\Bbb P$. Let $\overline{p}\in
\bigcap_{n\in\omega}D_{n}$ such that $\overline{p}\leq p$. Then
\[h=\{\langle n,a\rangle\in\omega\times A :
\overline{p}\forces\dot{h}(n)=a\}\in M.\]
$\Box$

\begin{lemma}

Let $S\subseteq\lim (\omega_{1})$ and
$T$ be an $\omega_{1}$--tree which is $S$--properly pruned in
countable products. Then for any index set $I$, the poset
${\Bbb P}(T,I,\omega_{1})$ is $\omega_{1}$--Baire.

\end{lemma}

\noindent {\bf Proof:}\quad
For each $n\in\omega$,
let $D_{n}$ be a dense open subset of ${\Bbb P}(T,I,\omega_{1})$. Let
$p\in {\Bbb P}(T,I,\omega_{1})$. We now construct $p_{s}\in {\Bbb P}
(T,I,\omega_{1})$ for every $s\in 2^{<\omega}$ inductively on the
length of $s$ such that,

(1) $p_{0}\leq p$, 

(2) $s\subseteq t$ iff $p_{t}\leq p_{s}$,

(3) there is an increasing sequence $\langle\delta_{n}:n\in\omega\rangle$
of countable ordinals such that for every $s\in 2^{n}$, $p_{s}$ is uniform
at $\delta_{n}$.

(4) for every $s\in 2^{<\omega}$, $p_{s\hat{\;}\langle 0\rangle}$ and
$p_{s\hat{\;}\langle 1\rangle}$ are completely incompatible,

(5) for every $s\in 2^{n}$, $p_{s}\in D_{n}$.

\medskip

Assume that we have already had $p_{s}$ for every $s\in 2^{<n}$.

Let $s\in 2^{n-1}$ and $q^{s}\in D_{n}$ such that $q^{s}\leq p_{s}$.
Let $l=0,1$. 
By Lemma 5, there are $q^{s}_{l}<q^{s}$ such that $q^{s}_{0}$ and
$q^{s}_{1}$ are completely incompatible. Let \[\delta_{n}=
\bigcup\{ht_{T}(q^{s}_{l}(i)):i\in I,s\in 2^{n-1}, l=0,1\}+1.\]
$\delta_{n}$ is countable because the support of every $q^{s}_{l}$
is countable. Let $p_{s\hat{\;}\langle l\rangle}\in {\Bbb P}(T,I,\omega_{1})$
such that $p_{s\hat{\;}\langle l\rangle}\leq q^{s}_{l}$ and 
all $p_{s\hat{\;}\langle l\rangle}$ are uniform at $\delta_{n}$.
$p_{s\hat{\;}\langle l\rangle}\in D_{n}$ because $p_{s\hat{\;}\langle l\rangle}
\leq q^{s}$.

Let $I'=\bigcup_{s\in 2^{<\omega}}supt(p_{s})$. Then $I'$ is countable.
\[C\!\res\!I'=\{p_{s}\!\res\!I':s\in 2^{<\omega}\}\] is now a Cantor tree in
$T^{I'}$, which is separated by (4) and uniformly cofinal in 
$(T\!\res\!\delta)^{I'}$, where $\delta=\bigcup_{n\in\omega}\delta_{n}$.
Since $T$ is $S$--properly pruned in countable
products, there exists $f\in 2^{\omega}$ such that
for every $i\in I'$, $\bigcup_{n\in\omega}p_{f\res n}(i)\in T_{\delta}
\bigcup\{\emptyset\}$.

Let $p_{f}\in {\Bbb P}(T,I,\omega_{1})$ defined by letting
\[p_{f}\!\res\!I'=\langle\bigcup_{n\in\omega}p_{f\res n}(i):i\in I'\rangle\]
and \[p_{f}\!\res\!I\smallsetminus I'\equiv\emptyset.\] Then $p_{f}\in
{\Bbb P}(T,I,\omega_{1})$ and $p_{f}\leq p_{f\res n}$ for every $n\in\omega$.
So $p_{f}\leq p$ and $p_{f}\in\bigcap_{n\in\omega}D_{n}$. \quad $\Box$

\begin{theorem}

Assume the existence of a strongly inaccessible cardinal. It is consistent
with $C\!H$ plus $2^{\omega_{1}}>\omega_{2}$ that there exists a Jech--Kunen
tree and there are no Kurepa trees.

\end{theorem}

\noindent {\bf Proof:}\quad
Let $M$ be a model of $G\!C\!H$ plus that there is a strongly inaccessible
cardinal $\kappa$. In $M$, let $T$ be an $\omega_{1}$--tree which is
$\lim (\omega_{1})$--properly pruned in countable products and let
$\mu$ and $\lambda$ be two regular cardinals such that $\kappa
\leq\mu<\lambda$. Again in $M$ let ${\Bbb P}_{1}=Lv(\kappa,\omega_{1})$,
${\Bbb P}_{2}={\Bbb P}(T,\mu,\omega_{1})$ and ${\Bbb P}_{3}=
Fn(\lambda,2,\omega_{1})$. Let $G=G_{1}\times G_{2}\times G_{3}$
be a ${\Bbb P}_{1}\times {\Bbb P}_{2}\times {\Bbb P}_{3}$--generic filter
over $M$. We will show that $M[G]$ is a model of $C\!H$ plus 
$\lambda=2^{\omega_{1}}>\mu\geq\omega_{2}=\kappa$, in which there
are no Kurepa tree and $T$ is a Jech--Kunen tree with $\mu$--many branches.

\medskip

{\bf Claim 8.1}.\quad $M^{\omega}\bigcap M[G]\subseteq M$.

Proof of Claim 8.1: \quad
We first force with ${\Bbb P}_{2}$. 
By Lemma 6 and Lemma 7, ${\Bbb P}_{2}$ is $\omega_{1}$--Baire and forcing
with ${\Bbb P}_{2}$ will not add any new countable sequences. Hence
${\Bbb P}_{1}\times {\Bbb P}_{3}$ is still $\omega_{1}$--closed in
$M[G_{2}]$. Then forcing with ${\Bbb P}_{1}\times {\Bbb P}_{3}$ will
also not add any new countable sequences because it is ${\omega_{1}}$--closed.

\medskip

{\bf Claim 8.2}.\quad ${\Bbb P}_{1}\times {\Bbb P}_{2}\times {\Bbb P}_{3}$
has the $\kappa$--c.c..

Proof of Claim 8.2:\quad Let 
\[\{\langle p_{\alpha},q_{\alpha},r_{\alpha}\rangle:
\alpha<\kappa\}\subseteq {\Bbb P}_{1}\times {\Bbb P}_{2}\times {\Bbb P}_{3}.\]
By the $\Delta$--system lemma, we can assume that the domains of all
$p_{\alpha}$, the domains of all $q_{\alpha}$ and the domains of all
$r_{\alpha}$ form three $\Delta$--systems with roots $\Delta_{1}$,
$\Delta_{2}$ and $\Delta_{3}$ respectively. Since there are less than
$\kappa$--many $p$'s in ${\Bbb P}_{1}$ with domains $=\Delta_{1}$,
there are $\omega_{1}$--many $q$'s in ${\Bbb P}_{2}$ 
with domains $=\Delta_{2}$, 
and there are $\omega_{1}$--many $r$'s in ${\Bbb P}_{3}$ with domains
$=\Delta_{3}$, then there exist $\alpha_{1}$ and $\alpha_{2}$
in $\kappa$ such that 
\[p_{\alpha_{1}}\!\res\!\Delta_{1}=
p_{\alpha_{2}}\!\res\!\Delta_{1},\;
q_{\alpha_{1}}\!\res\!\Delta_{2}=q_{\alpha_{2}}\!\res\!
\Delta_{2}\mbox{ and }r_{\alpha_{1}}\!\res\!\Delta_{3}=
r_{\alpha_{2}}\!\res\!\Delta_{3}.\]
Obviously $\langle p_{\alpha_{1}},q_{\alpha_{1}},r_{\alpha_{1}}\rangle$
and $\langle p_{\alpha_{2}},q_{\alpha_{2}},r_{\alpha_{2}}\rangle$ are
compatible.

\medskip

\noindent Remark: \quad By Claim 8.1 and Claim 8.2, $\omega_{1}$ and 
all the cardinals 
greater than or equal to $\kappa$ in $M$ are preserved and 
$C\!H$ is true in $M[G]$. 
In $M[G]$, $\kappa=\omega_{2}$ because forcing with ${\Bbb P}_{1}$
collepses all the cardinals between $\omega_{1}$ and $\kappa$ in $M$.
Also in $M[G]$, $2^{\omega_{1}}=\lambda$ because forcing with ${\Bbb P}_{3}$
adds $\lambda$--many subsets of $\omega_{1}$. 

\medskip

{\bf Claim 8.3}.\quad There are no Kurepa trees in $M[G]$.

Proof of Claim 8.3: \quad Suppose that is not true. Let $K$ be a normal
Kurepa tree in $M[G]$. Since $|K|=\omega_{1}$, there are $\theta<\kappa$,
$I\subseteq\mu$ with $|I|\leq\omega_{1}$ and $J\subseteq\lambda$ with
$|J|\leq\omega_{1}$ such that 
\[K\in M[G']=M[G'_{1}\times G'_{2}\times
G'_{3}],\] where \[G'_{1}=G_{1}\bigcap Lv(\theta,\omega_{1}),\]
\[G'_{2}=G_{2}\bigcap {\Bbb P}(T,I,\omega_{1}),\] \[G'_{3}=G_{3}\bigcap
Fn(J,2,\omega_{1})\] and 
\[G'=G'_{1}\times G'_{2}\times G'_{3}.\]
Let \[G''_{1}=G_{1}\bigcap Lv(\kappa\smallsetminus\theta,
\omega_{1}),\]
\[G''_{2}=G_{2}\bigcap 
{\Bbb P}(T,\mu\smallsetminus I,\omega_{1}),\]
\[G''_{3}=G_{3}\bigcap Fn(\lambda\smallsetminus J,2,\omega_{1})\] and
\[G''=G''_{1}\times G''_{2}\times G''_{3}.\]
Since $M[G']\models 2^{\omega_{1}}<\kappa$, there exists \[b\in
{\cal B}(K)\bigcap M[G]\smallsetminus M[G'].\] Furthermore \[b\not\in
M[G'][G''_{1}][G''_{3}]\]
because $Lv(\kappa\smallsetminus\theta,\omega_{1})$
and $Fn(\lambda\smallsetminus J,2,\omega_{1})$ are $\omega_{1}$--closed
in $M[G']$. We now work in $M[G'][G''_{1}][G''_{3}]$ and
let $p\in G''_{2}$ such that \[p\forces (\dot{b}\in {\cal B}(K)\bigcap
M[\dot{G}]\smallsetminus M[G'][G''_{1}][G''_{3}]).\]

\medskip

We construct \[C=\{p_{s}:s\in 2^{<\omega}\}\subseteq {\Bbb P}(T,
\mu\smallsetminus I,\omega_{1})\] and \[D=\{k_{s}:s\in 2^{<\omega}\}
\subseteq K\] such that,

(1) $s\subseteq s'$ iff $p_{s'}\leq p_{s}$ iff $k_{s}\leq k_{s'}$,

(2) $C$ is separated and uniformly cofinal in 
$(T\!\res\!\delta)^{\mu\smallsetminus I}$ for some $\delta\in\lim (\omega_{1})$,

(3) $D$ is cofinal in $K\!\res\!\delta'$ for some $\delta'\in\lim (\omega_{1})$,

(4) for every $s\in 2^{<\omega}$, $p_{s}\forces k_{s}\in\dot{b}$.

\medskip

Assume that we have already had $p_{s}$ and $k_{s}$ for all $s\in 2^{<n}$.
Let \[\delta'_{n}=\bigcup\{ht_{K}(k_{s}):s\in 2^{<n}\}+1\] and pick $s\in
2^{n-1}$. Let $l=0,1$.

First find $p'_{s}\leq p_{s}$ such that \[\exists x\in K_{\delta'_{n}}
(p'_{s}\forces x\in\dot{b}).\] Since \[p'_{s}\forces\dot{b}\not\in
M[G'][G''_{1}][G''_{3}],\] there exist $q^{s}_{l}\leq p'_{s}$ and 
$x_{l}>x>k_{s}$ such that $x_{0}\perp x_{1}$ and $q^{s}_{l}\forces
x_{l}\in\dot{b}$. By Lemma 5, we can extend $q^{s}_{l}$ to $r^{s}_{l}$
such that $r^{s}_{l}$ are uniform at $\alpha_{s}<\omega_{1}$ and
$r^{s}_{0}$ is completely incompatible with $r^{s}_{1}$.
Let \[\delta_{n}=\bigcup\{\alpha_{s}:s\in 2^{n-1}\}+1,\] 
$p_{s\hat{\;}\langle l\rangle}$ be an extension of $r^{s}_{l}$ such that
$supt(p_{s\hat{\;}\langle l\rangle})=supt(r^{s}_{l})$ and
$p_{s\hat{\;}\langle l\rangle}$ be uniform at $\delta_{n}$.
This ends the construction.

\medskip

Let $\delta=\bigcup_{n\in\omega}\delta_{n}$, $\delta'=\bigcup_{n\in\omega}
\delta'_{n}$ and
$I'=\bigcup_{s\in 2^{<\omega}}supt(p_{s})$. Then $I'$ is countable.
Since $T$ is $\lim (\omega_{1})$--properly pruned in countable products
and $C\!\res\!I'$ is a Cantor tree which is separated and uniformly cofinal
in $(T\!\res\!\delta)^{I'}$, then there are uncountably many $f\in 2^{\omega}$
such that $p_{f}$ defined by letting \[p_{f}(i)=\bigcup_{n\in\omega}
p_{f\res n}(i)\] for every $i\in I'$ is a lower bound of $\{p_{f\res n}:
n\in\omega\}$ in ${\Bbb P}(T,\mu\smallsetminus,\omega_{1})$.
(Note that $C$ is in $M$ because no new countable sequences are added.)
For every such $f$ there exists $k_{f}\in K_{\delta'}$ such that
$p_{f}\forces k_{f}\in\dot{b}$ 
and for different $f$, $k_{f}$ are different. That contradicts that $K$ is
a Kurepa tree.

\medskip

{\bf Claim 8.4}.\quad $M[G]\models (|{\cal B}(T)|=\mu)$.

Proof of Claim 8.4:\quad
$|{\cal B}(T)|\geq\mu$ is trivial because forcing with ${\Bbb P}_{2}$
adds at least $\mu$--many new branches of $T$. 
Since in $M[G_{1}][G_{2}]$, $2^{\omega_{1}}=\mu$, then we need only to show
that forcing with ${\Bbb P}_{3}$ will not add any new branches of $T$.

Suppose that is not true and let $b$ be a branch of $T$, which is
in $M[G]\smallsetminus M[G_{1}][G_{2}]$. 
We now work in $M[G_{1}][G_{2}]$ and let $p\in G_{3}$ such that
\[p\forces \dot{b}\in {\cal B}(T)\bigcap 
M[\dot{G}]\smallsetminus M[G_{1}][G_{2}].\]
We can then easily construct $C=\{p_{s}:s\in 2^{<\omega}\}\subseteq 
{\Bbb P}_{3}$ and $D=\{t_{s}:s\in 2^{<\omega}\}\subseteq T$ such that,

(1) $s\subseteq s'$ iff $p_{s'}\leq p_{s}$ iff $t_{s}\leq t_{s'}$,

(2) $D$ is a Cantor tree which is cofinal in $T\!\res\!\delta$ for some
$\delta\in\lim (\omega_{1})$,

(3) for every $s\in 2^{<\omega}$, $p_{s}\forces t_{s}\in\dot{b}$.

\medskip

Since $T$ is $\lim (\omega_{1})$--properly pruned by Lemma 3, there exists
$g\in  2^{\omega}$ such that $\bigcup_{n\in\omega}t_{g\res n}\not\in
T_{\delta}$. But ${\Bbb P}_{3}$ is $\omega_{1}$--closed
in $M[G_{1}][G_{2}]$ because no new countable sequences have been added.
Hence there exists $p_{f}\in {\Bbb P}_{3}$ such that
$p_{g}\leq p_{g\res n}$ for every $n\in\omega$. This implies that
there exists $t\in T_{\delta}$ such that $p_{f}\forces t\in\dot{b}$.
Hence \[t=\bigcup_{n\in\omega}t_{g\res n}\in T_{\delta},\] a contradiction.
\quad $\Box$

\bigskip

In the model constructed above, there is only one Jech--Kunen tree. Next we will
build a model of $C\!H$ plus $2^{\omega_{1}}>\omega_{2}$, in which there
are no Kurepa trees and there are many Jech--Kunen trees with different
numbers of branches.

\begin{theorem}

Assume the existence of a strongly inaccessible cardinal. It is consistent
with $C\!H$ plus $2^{\omega_{1}}>\omega_{2}$ that there are no Kurepa trees
and there are Jech--Kunen trees $T^{\alpha}$ for $\alpha\in\omega_{1}$
such that $\alpha\neq\alpha'$ implies $|{\cal B}(T^{\alpha})|\neq
|{\cal B}(T^{\alpha'})|$.

\end{theorem} 

\noindent {\bf Proof:}\quad
Let $M$ be a model of $G\!C\!H$ and that there exists a strongly inaccessible
cardinal $\kappa$. In $M$, let \[\Gamma=\{\mu_{\alpha}:\alpha\in\omega_{1}\}
\subseteq [\kappa,\lambda)\] be a set of different regular cardinals,
where $\lambda$ is also a regular cardinal. Again in $M$, let $\{S_{\alpha}:
\alpha\in\omega_{1}\}$ be a partition
of $\lim (\omega_{1})$ such that every $S_{\alpha}$ is a stationary, 
and let $T^{\alpha}$ be an $\omega_{1}$--tree which is
$S_{\alpha}$--properly pruned in countable products for every $\alpha\in
\omega_{1}$. In $M$, let ${\Bbb P}_{1}=Lv(\kappa,\omega_{1})$,
${\Bbb P}_{2}$ be the product of $\{{\Bbb P}(T^{\alpha},\mu_{\alpha},
\omega_{1}):\alpha\in\omega_{1}\}$ with countable supports, and
${\Bbb P}_{3}=Fn(\lambda,2,\omega_{1})$.
Let $G=G_{1}\times G_{2}\times G_{3}$ be a ${\Bbb P}_{1}\times {\Bbb P}_{2}
\times {\Bbb P}_{3}$--generic filter over $M$. Then $M[G]$ is the model 
we are looking for.

\medskip

{\bf Claim 9.1}.\quad $M^{\omega}\bigcap M[G]\subseteq M$.

\medskip

{\bf Claim 9.2}. \quad ${\Bbb P}_{1}\times {\Bbb P}_{2}\times {\Bbb P}_{3}$
has the $\kappa$--c.c..

\medskip

{\bf Claim 9.3}.\quad There are no Kurepa trees in $M[G]$.

\medskip

All the proofs of above three claims are similar to the proofs of
corresponding claims in Theorem 8. By Claim 9.1 and Claim 9.2, 
$\omega_{1}$ and all the cardinals greater than or equal to $\kappa$ are
preserved. Besides, forcing with ${\Bbb P}_{1}$ collapses all the cardinals
between $\omega_{1}$ and $\kappa$. So in $M[G]$, $C\!H$ is true,
$\kappa=\omega_{2}<\lambda =2^{\omega_{1}}$
and $\{\mu_{\alpha}:\alpha\in\omega_{1}\}\subseteq [\kappa,\lambda)$ is
still a set of different cardinals.

\medskip

{\bf Claim 9.4}.\quad $M[G]\models (|{\cal B}(T^{\alpha})|=\mu_{\alpha})$
for every $\alpha\in\omega_{1}$.

Proof of Claim 9.4: \quad Pick an $\alpha\in\omega_{1}$.
Let ${\Bbb P}^{\alpha}_{2}={\Bbb P}(T^{\alpha},\mu_{\alpha},\omega_{1})$
and ${\Bbb P}^{-\alpha}_{2}$ be the product of $\{{\Bbb P}(T^{\beta},
\mu_{\beta},\omega_{1}):\beta\in\omega_{1}\smallsetminus\{\alpha\}\}$
with countable supports.
Then ${\Bbb P}_{2}\cong {\Bbb P}^{\alpha}_{2}\times {\Bbb P}^{-\alpha}_{2}$.
Let $p\in {\Bbb P}_{2}$. We let 
\[S\!U\!P\!T(p)=\{\alpha\in\omega_{1}:supt(p(\alpha))
\neq\emptyset\}.\] Notice the differences between $supt$ and $S\!U\!P\!T$.
We call an element $p\in {\Bbb P}^{-\alpha}_{2}$ uniform at $\gamma$ for some
$\gamma\in\omega_{1}$ if for every $\beta\in S\!U\!P\!T(p)$, $p(\beta)$ is
uniform at $\gamma$.

\medskip

{\bf Subclaim 9.4.1}.\quad Forcing with ${\Bbb P}^{-\alpha}_{2}$ will not add
any new branches to $T^{\alpha}$.

Proof of Subclaim 9.4.1:\quad Let $G_{2}=G^{\alpha}_{2}\times G^{-\alpha}_{2}
\subseteq {\Bbb P}^{\alpha}_{2}\times {\Bbb P}^{-\alpha}_{2}$.
Suppose that Subclaim 1 is not true and let $b$ be a branch of $T^{\alpha}$
shch that \[b\in M[G_{1}][G_{2}]\smallsetminus 
M[G_{1}][G^{\alpha}_{2}].\]

We now work in $M[G_{1}][G^{\alpha}_{2}]$ and
let $p\in G^{-\alpha}_{2}$ such that \[p\forces\dot{b}\in {\cal B}(T^{\alpha})
\smallsetminus M[G_{1}][G^{\alpha}_{2}].\] We construct recursively
a normal subtree $T'$ of $T^{\alpha}$ with every level countable, and
a subset $C=\{p_{t}:t\in T'\}$ of ${\Bbb P}^{-\alpha}_{2}$ such that,

(1) for every $\delta\in\omega_{1}$ there is $\gamma_{\delta}$ such that
$T'_{\delta}\subseteq T^{\alpha}_{\gamma_{\delta}}$,

(2) if $\delta\in\lim (\omega_{1})$, then
$\gamma_{\delta}=\bigcup_{\beta<\delta}
\gamma_{\beta}$,

(3) $p_{\emptyset}\leq p$, and for any $t,t'\in T'$, $t\leq t'$ iff
$p_{t'}\leq p_{t}$,

(4) for every $t\in T'_{\delta}$, there is $\gamma',\;\gamma_{\delta}\leq
\gamma'\leq\gamma_{\delta+1}$
such that $p_{t}$ is uniform at $\gamma'$,

(5) if $t\in T'_{\delta}$ for some $\delta\in\lim (\omega_{1})$,
then $p_{t}$ is uniform
at $ht_{T^{\alpha}}(t)=\gamma_{\delta}$,

(6) $t\perp t'$ implies that $p_{t}(\beta)$ and $p_{t'}(\beta)$
are completely incompatible for every $\beta\in S\!U\!P\!T(p_{t})\bigcap
S\!U\!P\!T(p_{t'})$,

(7) for every $t\in T'$, $p_{t}\forces t\in\dot{b}$.

\medskip

Assume that we have already had $T'\!\res\!\delta$ and $C\!\res\!\delta=
\{p_{t}:t\in T'\res\delta\}$.

\medskip

Case 1. \quad $\delta=\gamma+1$ for some $\gamma\in\omega_{1}$.

Pick $t\in T'_{\gamma}$ and let $l=0,1$.

Since \[p_{t}\forces\dot{b}\not\in M[G_{1}][G^{\alpha}_{2}],\] there
exist $t_{l}\in T^{\alpha}, t_{l}>t$ and $q^{t}_{l}\leq p_{t}$ such that
\[t_{0}\perp t_{1}\mbox{ and }q^{t}_{l}\forces t_{l}\in\dot{b}.\] Without
loss of generality we can pick $t_{l}$ such that $ht_{T^{\alpha}}(t_{l})
=\delta'$ for every $t\in T'_{\gamma}$ and $l=0,1$, where \[\delta'
>\bigcup\{\gamma'':p_{t}\mbox{ is uniform at }\gamma''\mbox{ for some }
t\in T'_{\gamma}\}.\] Besides, we can require that $q^{t}_{0}(\beta)$
and $q^{t}_{1}(\beta)$ are uniform and 
are completely incompatible for every $t\in
T'_{\gamma}$ and $\beta\in S\!U\!P\!T(q^{t}_{0})\bigcap S\!U\!P\!T(q^{t}_{1})$.
Let $\gamma'\in\omega_{1}$ such that $\gamma'>\delta'$ and \[\gamma'
>\bigcup\{\gamma'':q^{t}_{l}\mbox{ is uniform at }\gamma''\mbox{ for some }
t\in T'_{\gamma}\mbox{ and }l=0,1\}.\] 
Let $T'_{\delta}=\{t_{l}:t\in T'_{\gamma},l=0,1\}$ 
and let $p_{t_{l}}\leq q^{t}_{l}$ such
that $p_{t_{l}}$ is uniform at $\gamma'$.

\medskip

Case 2.\quad $\delta\in\lim (\omega_{1})$.

First $\gamma_{\delta}$ can't be in $S_{\alpha}$ because otherwise
every $T^{\beta}$ for $\beta\in\omega_{1}\smallsetminus\{\alpha\}$
is complete at level $\gamma_{\delta}$. But 
in $M[G_{1}][G^{\alpha}_{2}]$, $T^{\alpha}$ is still properly
pruned at level $\gamma_{\delta}$
because forcing with ${\Bbb P}_{1}\times {\Bbb P}_{2}^{\alpha}$ adds no 
new countable sequences, so that there exists 
$B\in {\cal B}(T'\!\res\!\delta)$ such that $B$ has no upper bound in
$T^{\alpha}$. On the other hand, $\{p_{t}:t\in B\}$ has a lower bound
$p_{B}$ in ${\Bbb P}^{-\alpha}_{2}$. Then \[p_{B}\forces\exists t\in
T^{\alpha}_{\gamma_{\delta}}(t\in\dot{b})\] implies that $B$ has an upper bound
in $T^{\alpha}$, a contradiction.

Assume that $\gamma_{\delta}\in S_{\beta}$ for some $\beta\neq\alpha$.
Since in $M[G_{1}][G^{\alpha}_{2}]$,
$T^{\beta}$ is properly pruned at level $\gamma_{\delta}$, then
for every $t\in T'\!\res\!\delta$ there exists 
$B_{t}\in {\cal B}(T'\!\res\!\delta)$
such that $t\in B_{t}$ and $\langle\bigcup_{t'\in B_{t}}p_{t'}(\beta)(i)
\rangle_{i\in\mu_{\beta}}\in {\Bbb P}(T^{\beta},\mu_{\beta},\omega_{1})$.

Now every $T^{\beta'}$ is complete at level $\gamma_{\delta}$ for $\beta'
\neq\beta$. We can define $p_{B_{t}}\in {\Bbb P}^{-\alpha}_{2}$ by letting
\[p_{B_{t}}(\beta)(i)=\bigcup_{t'\in B_{t}}p_{t'}(\beta)(i)\] for every
$\beta\in\omega_{1}\smallsetminus\{\alpha\}$ and $i\in\mu_{\beta}$.

Let $T'_{\delta}=\{\bigcup B_{t}:t\in T'\!\res\!\delta\}$ and let
$p_{\bigcup B_{t}}=p_{B_{t}}$. This ends the construction.

\medskip

Since $S_{\alpha}$ is stationary and by (2), $\{\gamma_{\delta}:
\delta\in\omega_{1}\}$ is a club set, then there exists 
$\delta\in\omega_{1}$ such that $\gamma_{\delta}\in S_{\alpha}$.
But this has been shown impossible.

\medskip

{\bf Subclaim 9.4.2}.\quad Forcing with ${\Bbb P}_{3}$ will not add any
new branches to $T^{\alpha}$.

Proof of Subclaim 9.4.2:\quad Similar (but much easier) to the proof of
Subclaim 9.4.1.

\medskip

By Subclaim 9.4.1 and Subclaim 9.4.2, all the branches of $T^{\alpha}$ in $M[G]$
are already in $M[G_{1}][G^{\alpha}_{2}]$. But in $M[G_{1}][G^{\alpha}_{2}]$
$2^{\omega_{1}}=\mu_{\alpha}$. So $|{\cal B}(T^{\alpha})|=\mu_{\alpha}$.
\quad $\Box$

\bigskip

\noindent {\bf Concluding remarks}.\quad (1) 
$\mu,\mu_{\alpha}$ and $\lambda$ are
not necessarily regular. 

(2) In Theorem 9, we can also have larger number of trees. For this
we use $S_{\alpha}$'s which are only almost disjoint.

(3) In the proof of Theorem 9, if we do not want to use stationary
sets, we can force the trees as part of the forcing, and then
prove that they are ``pruned together'', so using the stationary sets
simplifies the matter.

(4) We have used ($\lim (\omega_{1})\smallsetminus S$)--complete tree $T$
({\it i. e.} every branch of $T\!\res\!\delta$ for $\delta\in\lim 
(\omega_{1})\smallsetminus S$ has an upper bound in $T$). Our consideration
leads naturally to $S$--Kurepa trees. $T$ is called an $S$--Kurepa tree
if $T=\bigcup_{\alpha\in\omega_{1}}T(\alpha)$, where $ht(T(\alpha))=\alpha$,
$T(\alpha)=\bigcup_{\beta<\alpha}T(\beta)$ if $\alpha\in\lim (\omega_{1})$
and $|T_{\alpha}\bigcap\{\bigcup B :B\in {\cal B}(T(\alpha))\}|\leq\omega$
if $\alpha\in S$. So we may well consider $S$--Kurepa and
$(\lim (\omega_{1})\smallsetminus S)$--complete trees.

(5) The $T$ we build are not only 
$(\lim (\omega_{1})\smallsetminus S)$--complete, but also strongly proper
(see [S1] or/and [S2]).

\bigskip

Institute of Mathematics, 

The Hebrew University, 

Jerusalem, Israel.

\bigskip

Department of Mathematics, 

Rutgers University, 

New Brunswick, NJ, 08903, USA.

\bigskip

Department of Mathematics,

University of Wisconsin, 

Madison, WI 53706, USA.

\bigskip

{\em Sorting:} The first two addresses are the first author's; the last
one is the second author's. 

\end{document}